%% file: contcar5.tex
\documentclass[twoside,openright,11pt,a4paper,epsf]{article}
\usepackage{latexsym}
\usepackage[T1]{fontenc}
\usepackage{amsmath}
\usepackage{amsthm}
\usepackage{amsfonts}
\usepackage{amssymb}
\usepackage{epsfig}
\usepackage[all]{xy}

\newcommand{\color}[6]{}

\newcommand{\R}{\mathbb{R}}

\newcommand{\C}{\mathbb{C}}

\renewcommand{\P}{\mathbb{P}}

\newcommand{\cc}{\mathbb{\mathcal C}}

\newcommand{\ce}{\mathcal E}

\newcommand{\nbd}{neighbourhood }
\newcommand{\nbds}{neighbourhoods }

\newcommand{\priv}{\backslash}
\newcommand{\lra}{\longrightarrow}

\newcommand{\om}{\omega}
\newcommand{\eps}{\varepsilon}
\renewcommand{\phi}{\varphi}

\newcommand{\wdt}[1]{\widetilde{#1}}
\newcommand{\cqfd}{\hfill $\square$ \vspace{0.1cm}\\ }

\newcommand{\im}{\textnormal{Im}\,}
\newcommand{\re}{\textnormal{Re}\,}

\newtheorem{definition}{Definition}[section]
\newtheorem{thm}{Theorem}
\newtheorem{prop}[definition]{Proposition}
\newtheorem{lemma}[definition]{Lemma}

\title{$\cc^0$-rigidity of characteristics in symplectic geometry.}
\author{Emmanuel Opshtein.}
\date{}
\begin{document}
\maketitle
\begin{abstract}
The paper concerns a $C^0$-rigidity result for the charcteristic foliations in symplectic geometry. A 
symplectic homeomorphism (in the sense of Eliashberg-Gromov) which preserves a smooth hypersurface also preserves its characteristic foliation.
\end{abstract}
\section*{Introduction}
Gromov and Eliashberg showed that a $\cc^0$-limit of symplectic diffeomorphisms which is itself a diffeomorphism 
is symplectic (\cite{gromov,eliashberg}, see also \cite{hoze}). This rigidity result leads to the definition of symplectic homeomorphisms 
(the $\cc^0$-limits of symplectic diffeomorphisms which are homeomorphisms), and shows that they define a 
proper subset of volume preserving homeomorphisms in dimension greater than $4$. It also raises the question of the survival of the symplectic invariants to this limit process. Which classical invariants of symplectic geometry remain 
invariants of this maybe softer $\cc^0$-symplectic geometry?  This paper shows that the characteristic foliation is one of them. 
\begin{thm}\label{rigidity}
Let $S$ and $S'$ be smooth hypersurfaces of some symplectic manifolds $(M,\om)$, $(M',\om')$.  Any symplectic homeomorphism between $M$ and $M'$ which sends $S$ to $S'$ transports the characteristic foliation of $S$ to that of $S'$.
\end{thm}
The characteristic foliation is a symplectic invariant of a given hypersurface $S$, which can be defined as the integral foliation of the (one dimensional) null space of the restriction of the symplectic form to  $S$. This definition
is intrinsically smooth since it involves the tangent spaces of $S$. But the roles of this foliation in symplectic geometry are many. In particular, one of its rather folkloric property concerns non-removable intersection : 
if two smoothly bounded open sets intersect exactly on their boundaries, and if no symplectic perturbation can separate them, then the boundaries share a common closed invariant subset of the characteristic foliation 
\cite{lamc,lasi,moi3,popasi}.  This paper proceeds from the remark that this property have a meaning also in the continuous 
category, so defining this foliation in continuous terms is conceivable. 

An application of this theorem is a weak answer to a question by Eliashberg and Hofer  about the symplectic characterization of a hypersurface by the open set it bounds : {\it  Under which conditions the existence of a symplectomorphism between two smoothly bounded open sets in symplectic manifolds imply that their boundaries 
are symplectomorphic also \cite{elho} ?}  Some results are known \cite{cieliebak,buhi}, but nothing when the sets are standard balls. Theorem \ref{rigidity} allows a  partial answer in this case.
\begin{thm}\label{weakelho}
Let $U$ be a smoothly bounded open set in $\R^4$. Assume that there is a symplectomorphism between $B^4(1)$
and $U$ which extends continuously to a homeomorphism between $S^3$ and $\partial U$. Then $\partial U$
is symplectomorphic to $S^3$.
\end{thm}

The paper is organized as follows. We first define symplectic hammers and explain their roles : theorem \ref{rigidity} proceeds from a localization of their actions along characteristics (section \ref{hammers}).  
This localization is proved in section  \ref{theproof}. We then present the application in the last section. \vspace{,2cm}\\
{\bf Aknowledgements.} I wish to thank Leonid Polterovich for making me aware of a serious mistake in the first 
version of the proof.
\section{Symplectic hammers.}\label{hammers}
Let $S$ be a hypersurface in a symplectic manifold $M$. We say that   $B$ is a small ball centered on $S$ if   it is 
a symplectic embedding of an euclidean ball centered at the origin into $M$ which sends $\R^{2n-1}:=\R\times \C^{n-1}\subset \C^n$ to $S$.  Such a ball is disconnected by $S$ into two components denoted by $S_+$ and $S_-$. 
By a classical result, any point of $S$ is the center of such a ball. Fix also a metric on $M$ in order to refer to small  sets.
\begin{definition}\label{defhammer}
Given two points $x$, $y$ on $S\cap B$ and a (small) positive  real $\eps$, an  $\eps$-symplectic hammer between $x$ and $y$ with support in $B$ is a continuous path  of symplectic homeomorphisms $\Phi_t$ ($t\in[0,1]$) with common supports in $B$, and for which there exist two open sets $U_\eps(x)$ and $U_\eps(y)$ contained in the $\eps$-balls around $x$ and $y$ respectively 
such that : 
\begin{enumerate}
\item $\Phi_0=\,$Id,
\item  $\Phi_t(z)\in S_+$ for all $t\in]0,1]$ and $z\in S\cap U_\eps(x)$,
\item $\Phi_t(z)\in S_-$ for all $t\in]0,1]$ and  $z\in S\cap U_\eps(y)$,
\item $\Phi_t(z)\in S$ for all $t\in[0,1]$ and  $z\in S\priv \big(U_\eps(x)\cup U_\eps(y)\big)$.
\end{enumerate}
A smooth hammer will refer to a smooth isotopy of smooth symplectomorphisms verifying the four conditions above.
\end{definition}
\noindent In other terms, $\Phi_t$ preserves the hypersurface $S$ except for two bumps in opposite sides (a symmetry is necessary in view of the volume preservation). 

\begin{figure}[h]
\begin{center}
\input hamilmer.pstex_t
\end{center}
\end{figure}

\begin{figure}[h]
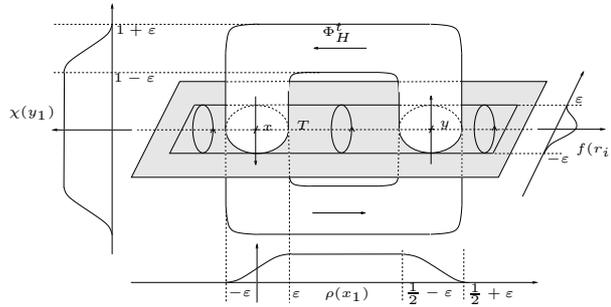

\begin{center}
\input example2.pstex_t
\caption{The Hamiltonian flow of $H$ in the proof of proposition \ref{existence}.}\label{example}
\end{center}
\end{figure}

\noindent One can easily construct examples of  symplectic hammers.
\begin{prop}\label{existence}
If $x,y\in B\cap S$ lie in the same characteristic, there exist $\eps$-symplectic hammers between $x$ and $y$ for all $\eps>0$.
\end{prop}
\noindent{\it Proof :} Since all hypersurfaces are locally symplectically the same, it is enough to produce a symplectic hammer for $\R^{2n-1}=\{\im z_1=0\}\subset \C^n$ between the points $p=0$ and $q=(1/2,0,\dots,0)$. Putting $x_1=\re z_1$, $y_1=\im z_1$ and $r_i=|z_i|$, consider a  Hamiltonian of the following type.
 $$
 H(z_1,\dots,z_n):= \chi(y_1)\rho(x_1)\Pi_{i=2}^n f(r_i).
 $$
If $\chi$, $\rho$ and $f$ are the bell functions represented in figure \ref{example},  and maybe multiplying $H$ by a small constant in order to slow the flow down produces  a symplectic hammer between $x$ and $y$.\hfill $\square$

 Proposition \ref{existence} can easily be reversed in the smooth category : two points lie in the same characteristic leaf of $S\cap B$ if and only if there exist smooth symplectic hammers between them. It is less obvious, but still true  that {\it all} the symplectic hammers also meet this constraint. Theorem \ref{rigidity} obviously follows because the class of symplectic hammers is preserved by symplectic homeomorphisms. 
\begin{prop}\label{c0carcar}
A hypersurface $S$ and a small ball $B$ centered on $S$ being given, there exists an $\eps$-hammer between $x,y \in B\cap S$ with support in $B$ for all small $\eps$ if and only if $x$ and $y$ are on the same characteristic leaf.
\end{prop}
\noindent{\it Proof of theorem \ref{rigidity} (assuming proposition \ref{c0carcar}) :} 
Let $M,M',S,S'$ and $\Phi$ be as in theorem \ref{rigidity}, and put any metric on $M$ and $M'$. 
Consider two points $x,y\in S$ which lie in the same  characteristic. Consider a covering $\cal B=\{B_\alpha\}$ of 
$S$ by small balls (in the above sense) whose images by $\Phi$ are contained 
in small balls $B_\alpha'$ centered on $S'$. 
Let $(x_i)_{i\leq N}$ be a chain between $x$ and $y$ (that is $x_0=x$, $x_N=y$) such that 
$x_i$ and $x_{i+1}$ are always in a same ball $B_i$. Then there exist $\eps$-hammers $\Phi^{(\eps)}_t$ 
with supports in $B_i$ between $x_i$ and $x_{i+1}$ for all $\eps$. The isotopies $\Phi\circ\Phi^{(\eps)}_t\circ\Phi^{-1}$ define continuous $\delta(\eps)$-symplectic hammers with support in $B_i'$ between $\Phi(x_i)$ and $\Phi(x_{i+1})$, where 
$\delta(\eps)$ goes to zero with $\eps$. Therefore by proposition \ref{c0carcar}, $\Phi(x_i)$
 and $\Phi(x_{i+1})$ are on the same characteristic, so $\Phi(x)$ and $\Phi(y)$ are also on the same characteristic.
 \hfill $\square$

\section{Proof of proposition \ref{c0carcar}}\label{theproof}
The idea is the following. Since preserving a foliation is a local property, and since all hypersurfaces are locally the same in the symplectic world, we could translate the non-preservation of {\it one} characteristic by a symplectic homeomorphism to  the existence of a local, hence universal object (a hammer between  points on distinct characteristics) which would exist on all hypersurfaces. These continuous hammers  would allow to break 
intersections between open sets  as long as these intersections only 
consist of one characteristic. But some such intersections are known to be non-removable : the most famous one being
the intersection between the complement of the cylinder $Z(1)$ and the closed ball $B^{2n}(1)$.  

\begin{lemma}\label{spherenough}
If proposition \ref{c0carcar} does not hold, then for any point $x$ of the euclidean sphere $S^{2n-1}\subset \C^n$
and for any positive $\eps$, there exists a continuous $\eps$-symplectic hammer between $x$ and a point $y$ which does not lie in the characteristic circle passing through $x$.
\end{lemma}
\noindent {\it Proof :} Assume that proposition \ref{c0carcar} does not hold. Then there exists a small ball 
centered on a hypersurface $S$, two points $p,q\in S\cap B$  which are not in the same characteristic 
of $S\cap B$ and a 
family $\Phi^{(\eps)}:=(\Phi_t^{(\eps)})_{t\in[0,1]}$ of $\eps$-symplectic hammers with supports in $B$ 
between $p$ and $q$.
By definition of a small ball, there is a symplectic diffeomorphism $\Psi_1$ between $B$ and an euclidean ball 
$B_1\subset \C^n$ around the origin with $\Psi(S\cap B)=\R^{2n-1}\cap B_1$.
Then  $\Psi_1$ takes $\Phi{(\eps)}$ to an $\eps$-hammer between $\Psi_1(p)$ and 
$\Psi_1(q)$ which are not on the same characteristic. By use of  translation and rescaling, we can assume 
that $\Psi_1(p)$ is the origin and $B_1$ is as small a \nbd of $0$ as wished.

Now given the point $x\in S^{2n-1}$, and if $B_1$ is small enough, there exists a symplectic diffeomorphism 
$\Psi_2:B_1\lra \C^n$ with $\Psi_2(B_1\cap \R^{2n-1})\subset S^{2n-1}$, $\Psi_2(0)=x$ and such that different characteristics
of $\R^{2n-1}\cap B_1$ are sent by $\Psi_2$ not only to different characteristics of $S^{2n-1}\cap \Psi_2(B_1)$ but even of $S^{2n-1}$ (this means that we do not allow $\Psi_2$ to "bend" $B_1$ so as to take two different characteristics to two different segments of a same characteristic circle of $S^{2n-1}$).  The continuous symplectic isotopies 
obtained by transporting $\Phi^{(\eps)}$ by $\Psi_2\circ \Psi_1$ are $\eps$-hammers between $x$ 
and the point $y:=\Psi_2(\Psi_1(q))$ which is not on the characteristic through $x$.\hfill $\square$

\begin{lemma}\label{sullivan}
Any bounded open set $U\subset Z(1):=B^2(1)\times \C^{n-1}$ whose boundary $\partial U$ does not contain a 
characteristic circle $S^1\times\{\cdot\}$ of $\partial Z(1)$ can be symplectically displaced from $\partial Z(1)$ 
to the interior of $Z(1)$.
\end{lemma}
\noindent{\it Proof :} Recall that the characteristic flow of $\partial Z(1)$ can be oriented by the vector field $JN$ where 
 the vector $N$ is the outward normal vector field to $\partial Z(1)$ and $J$ is the standard complex structure on $\C^n$. Now observe that if the compact set $K:=\partial U\cap \partial Z(1)$ does not contain any characteristic circle,
there exists a smooth function $H$ on $\C^n$ which decreases along the characteristic flow on a \nbd of $K$  \cite{sullivan}. The corresponding Hamiltonian vector field points toward the inside of the cylinder on this \nbd of $K$ because 
$$
g(X_H(x),N(x))=\om(X_H(x),JN(x))=dH(JN(x))<0,
$$
so $U$ is driven inside $Z(1)$ by the flow of $H$ for small enough times. 
\cqfd
\indent Let us come back to proposition \ref{c0carcar}. Consider $B^{2n}(1)$ as an open set lying in $Z(1)$. Its boundary $S^{2n-1}$ meets $Z(1)$ along precisely 
one characteristic circle of $\partial Z(1)$ :
$$
S^{2n-1}\cap \partial Z(1)=\{|z_1|=1,\; z'=0\}\subset B^2_{z_1}(1)\times \C^{n-1}_{z'}.
$$
Assume then by contradiction that proposition \ref{c0carcar} does not hold. 
Then by lemma \ref{spherenough}, there exists 
 $\eps$-hammers $\Phi_t$ between the point $(1,0)\in S^{2n-1}\cap \partial Z(1)$ and an interior  
point $y\in S^{2n-1}\cap Z(1)$ for arbitrarily small $\eps$. If $\eps$ is small enough, and slowing the flow of the hammer down enough (considering $\Phi_{at}$ in place of $\Phi_t$), 
the image $U$ of $B(1)$ by $\Phi_1$ is therefore an open set of $Z(1)$ whose boundary intersection with 
$\partial Z(1)$ is the circle $S^1\times \{0\}$ minus a small \nbd of $(1,0)$ which was taken inside the ball - hence inside the cylinder - by our hammer. By lemma \ref{sullivan}, there exists therefore a smooth Hamiltonian $K$ such that 
$\Phi^1_K(U)$ is relatively compact in $Z(1)$. Taking a good enough smooth approximation 
$\wdt \Phi_t$ of $\Phi_t$, $\Phi^1_K\circ \wdt\Phi_1(B(1))$ is still relatively compact in $Z(1)$. But this is in contradiction with 
Gromov's non-squeezing theorem.\hfill $\square$  

\section{Symplectic geometry from the inside.}
In this section, we prove theorem \ref{weakelho}. We first prove that the characteristics of the sphere are sent to the characteristics of the boundary of $U$. The next point is to see that the action on these characteristics coincide on both hypersurfaces. The result follows.
\subsection{A one-sided version of theorem \ref{rigidity}.}
In fact, theorem \ref{rigidity} holds also in a slightly more general framework.
\begin{thm}\label{onesided}
Let $U$ and $U'$ be smoothly bounded open sets in symplectic manifolds. Any symplectic homeomorphism between 
$U$ and $U'$ which extends continuously to a homeomorphism of their closures transports the characteristic foliation
of $\partial U$ to that of $\partial U'$.
\end{thm} 
The proof below is rather quick because everything has already been explained. Exactly as for theorem \ref{rigidity},
the point is to define a convenient notion of symplectic hammer which is invariant by one-sided symplectic homeomorphisms. Note that definition \ref{defhammer}  {\it has} to be modified 
since it  involves both sides of the hypersurface through the bumps. The solution is simply to forget about the part of the hammer which goes outside $U$.
\begin{definition}[One-sided hammers]\label{oshammer} Let $U$ be a smoothly bounded open set in $M$, $B$ a small 
ball centered in $\partial U$, $x,y\in B\cap \partial U$. A one-sided $\eps$-symplectic hammer between $x$ 
and $y$ is a continuous isotopy of homeomorphism $\Phi_t:\overline{U}\priv B_\eps(y)\lra \overline{U}$ which can be uniformly approximated in $U\priv B_\eps(y)$ by smooth symplectic isotopies with common supports in $B$, and   which verifies also properties 1), 2) and 4) of definition \ref{defhammer}.  
\end{definition}
\noindent This definition actually provides a one-sided definition of the characteristics because of the following. 
\begin{prop}
There exists one-sided $\eps$-hammers between $x$ and $y$ for all small $\eps$ if and only if $x$ and $y$ belong 
to the same characteristic of $\partial U$.
\end{prop}  
\noindent {\it Proof :} The proof is very similar to the proof of proposition \ref{c0carcar}. On one hand, since the
definition of a one-sided hammer is a local one, and since smooth hypersurfaces have no local invariants, the 
existence of a one-sided hammer between two points not in the same characteristic of $\partial U\cap B$ ensures
the existence of such a hammer on the ellipsoid 
$$
\ce(1,2):=\{|z_1|^2+\frac{|z'|^2}{4}\leq 1\}\subset Z(1)\subset \C_z\times \C_{z'}
$$ 
between the point $(1,0)$ of the "least action" characteristic $C_0$ and another point $y$ not in this characteristic. This hammer isotops $\ce(1,2)\priv B_\eps(y)$ to an open set $U\subset Z(1)$ whose boundary contains no characteristic circle of $\partial Z(1)$. By a Hamiltonian flow, $U$ can be symplectically displaced from $\partial Z(1)$
inside $Z(1)$. For $\eps$ small enough, and since $y$ does not belong to $C_0$, $\ce(1,2)\priv B_\eps(y)$
contains the ball of radius $1$, contradicting again Gromov's non-squeezing theorem.\hfill $\square$

\subsection{Proof of theorem \ref{weakelho}.}
Let us fix the notations. On $S^3=\partial B^4(1)$, the characteristic foliation defines the Hopf fibration
$\pi:S^3\lra \P^1$. Moreover, if $\om_0$ is the standard area form on $\P^1$ with total 
area $\pi$, the restriction of the symplectic form $\om$ on $\R^4$ to $S^3$ is $\pi^*\om_0$.

Let $U$ be a smoothly bounded domain in $\R^4$,  $f$ a symplectic diffeomorphism from $U$ to $B^4(1)$ which extends continuously to a homeomorphism of the boundaries. Then by theorem \ref{onesided}, $f$ sends the characteristics of $\partial U$ to those of $S^3$, so the characteristic foliation of $\partial U$ is a topological 
Hopf fibration.  By a work of Epstein \cite{epstein}, there is a {\it diffeomorphism} 
$\Psi:\partial U\lra S^3$ which takes the characteristics of $\partial U$ to the Hopf circles of $S^3$. Since 
the restriction of $\Psi_*\om$ to $S^3$ vanishes along the Hopf circles but never vanishes, 
$\Psi_*\om$ is also the pull-back of an area form $\om_1$ on $\P^1$ : $\Psi_*\om=\pi^*\om_1$. We claim 
that the $\om_1$-area of $\P^1$ is the expected one :
\begin{equation}\label{final}
\left\vert \int_{\P^1}\om_1\right \vert=\pi.
\end{equation}
Let us accept this as a fact for a moment. Then, after a possible orientation-reversing change of coordinates on 
$\P^1$, there exists an area-preserving diffeomorphism 
$\phi:(\P^1,\om_1)\lra (\P^1,\om_0)$. Any lift $\Phi$ of $\phi$ through $\pi$ is a self-diffeomorphism of 
$S^3$ which pulls back $\pi^*\om_0=\om_{|S^3}$ to $\Psi_*\om$. 
$$
\xymatrix{
      (\partial U,\om) \ar[r]^\Psi & (S^3,\Psi_*\om)\ar[r]^\Phi \ar[d]^\pi & (S^3,\pi^*\om_0) \ar[d]^\pi\\
      &  (\P^1,\om_1)\ar[r]^\phi & (\P^1,\om_0)
    }
$$
Hence, $(\Phi\circ \Psi)^*\om_{|S^3}=\om_{|\partial U}$ and theorem \ref{weakelho}  
follows by a classical argument of standard \nbd (see \cite{mcsa}). 

In order to prove (\ref{final}), observe that putting 
$$
\left \vert\int_{\P^1} \om_1 \right \vert=\pi r^2,
$$
the above argument shows the existence of a symplectic diffeomorphism $\Phi$ between a \nbd of the euclidean sphere of radius $r$ and a \nbd of $\partial U$ which takes $S^3(r)$ to $\partial U$.  The 
map $g=f\circ\Phi$ defines therefore a symplectomorphism between two one-sided \nbds $V_r$, $V_1$ of the euclidean spheres. The following lemma finally ensures that $r=1$.\hfill $\square$
\begin{lemma}
There exists a symplectomorphism between two open one-sided \nbds of euclidean spheres  $S^3(r)$ and 
$S^3(R)$ if and only if $r=R$.
\end{lemma}
It is straightforward in view of \cite{elgr}. We give however  the argument for the convenience of the reader.   
Notice that it would be obvious, should the map $g$ extend smoothly to $S^3(r)$ (which is precisely not the case) because the action of the characteristics on $S^3(r)$ and $S^3$ should coincide. \\
{\it Proof :} Assume that $r\leq R$, call $V_r$ and $V_R$ the one-sided open \nbds of the two spheres and $g:V_R\lra V_r$ a symplectomorphism. Notice that since $V_R$ is symplectically convex with respect to  $S^3(R)$, so is $V_R$  with respect to $S^3(r)$ (meaning that there is a contracting vector field on $V_r$ flowing away from $S^3(r)$),  so both $V_r$ and $V_R$ are contained inside the corresponding euclidean balls. On $B^4(R)$, the contracting vector field $X=-\sum r_i\partial/\partial r_i$ is $\om$-dual to the form $\lambda_0=-\sum r_i^2d\theta_i$. 
Its image $g_*X$ is also a contracting vector field on $V_r$, $\om$-dual to $g_*\lambda_0$. Since $H^1(V_r)=0$, $g_*\lambda$
extends to a primitive of $-\om$ on $B^4(r)$.  This extension provides $B^4(r)$ with a contracting vector field $X'$ which coincides with $g_*X$ on $V_r$ and which is forward complete (its flow is defined for all positive time) because it points inside $B^4(r)$ near its boundary. Therefore $g$ can be extended to a symplectic embedding $\wdt g:B^4(R)\lra B^4(r)$ by the formula :
$$
\wdt g(p)=\Phi^t_{X'}\circ g \big(\Phi^{-t}_X(p)\big),  \; \forall p \in B^4(R), \; \forall t \text{ such that } \Phi^{-t}_X(p)\in V_R.
$$ 
But such an embedding is only possible if $R\leq r$ by volume considerations. \hfill $\square$

{\footnotesize
\bibliographystyle{abbrv}
\bibliography{bib3.bib}
}

\vspace{5cm}
\noindent Emmanuel Opshtein,\\
Université de Strasbourg,\\
IRMA, UMR 7501 \\
7 rue René-Descartes \\
67084 Strasbourg Cedex \\

\end{document}

%% file: hamilmer.pstex_t
\begin{picture}(0,0)%
\epsfig{file=hamilmer.pstex}%
\end{picture}%
\setlength{\unitlength}{1782sp}%
\begingroup\makeatletter\ifx\SetFigFont\undefined%
\gdef\SetFigFont#1#2#3#4#5{%
  \reset@font\fontsize{#1}{#2pt}%
  \fontfamily{#3}\fontseries{#4}\fontshape{#5}%
  \selectfont}%
\fi\endgroup%
\begin{picture}(8478,3233)(264,-5059)
\put(3715,-3375){\makebox(0,0)[lb]{\smash{{\SetFigFont{6}{7.2}{\rmdefault}{\mddefault}{\updefault}{\color[rgb]{0,0,0}$\Phi_t$}%
}}}}
\put(2562,-3704){\makebox(0,0)[lb]{\smash{{\SetFigFont{6}{7.2}{\rmdefault}{\mddefault}{\updefault}{\color[rgb]{0,0,0}$x$}%
}}}}
\put(5772,-3108){\makebox(0,0)[lb]{\smash{{\SetFigFont{6}{7.2}{\rmdefault}{\mddefault}{\updefault}{\color[rgb]{0,0,0}$y$}%
}}}}
\put(2853,-4992){\makebox(0,0)[lb]{\smash{{\SetFigFont{7}{8.4}{\rmdefault}{\mddefault}{\updefault}{\color[rgb]{0,0,0}$\Phi_t(x)$}%
}}}}
\end{picture}%

%% file: example2.pstex_t
\begin{picture}(0,0)%
\epsfig{file=example2.pstex}%
\end{picture}%
\setlength{\unitlength}{829sp}%
\begingroup\makeatletter\ifx\SetFigFont\undefined%
\gdef\SetFigFont#1#2#3#4#5{%
  \reset@font\fontsize{#1}{#2pt}%
  \fontfamily{#3}\fontseries{#4}\fontshape{#5}%
  \selectfont}%
\fi\endgroup%
\begin{picture}(17637,8932)(-3337,-7491)
\put(-3337,-1898){\makebox(0,0)[lb]{\smash{{\SetFigFont{5}{6.0}{\rmdefault}{\mddefault}{\updefault}{\color[rgb]{0,0,0}$\chi(y_1)$}%
}}}}
\put(13431,-3013){\makebox(0,0)[lb]{\smash{{\SetFigFont{5}{6.0}{\rmdefault}{\mddefault}{\updefault}{\color[rgb]{0,0,0}$f(r_i)$}%
}}}}
\put(6001,-7314){\makebox(0,0)[lb]{\smash{{\SetFigFont{5}{6.0}{\rmdefault}{\mddefault}{\updefault}{\color[rgb]{0,0,0}$\rho(x_1)$}%
}}}}
\put(5173,-2313){\makebox(0,0)[lb]{\smash{{\SetFigFont{5}{6.0}{\rmdefault}{\mddefault}{\updefault}{\color[rgb]{0,0,0}$T$}%
}}}}
\put(5941,479){\makebox(0,0)[lb]{\smash{{\SetFigFont{5}{6.0}{\rmdefault}{\mddefault}{\updefault}{\color[rgb]{0,0,0}$\Phi^t_H$}%
}}}}
\put(10226,-7314){\makebox(0,0)[lb]{\smash{{\SetFigFont{5}{6.0}{\rmdefault}{\mddefault}{\updefault}{\color[rgb]{0,0,0}$\frac{1}{2}+\eps$}%
}}}}
\put(-229,-909){\makebox(0,0)[lb]{\smash{{\SetFigFont{5}{6.0}{\rmdefault}{\mddefault}{\updefault}{\color[rgb]{0,0,0}$1-\eps$}%
}}}}
\put(-169,546){\makebox(0,0)[lb]{\smash{{\SetFigFont{5}{6.0}{\rmdefault}{\mddefault}{\updefault}{\color[rgb]{0,0,0}$1+\eps$}%
}}}}
\put(13406,-1509){\makebox(0,0)[lb]{\smash{{\SetFigFont{5}{6.0}{\rmdefault}{\mddefault}{\updefault}{\color[rgb]{0,0,0}$\eps$}%
}}}}
\put(12566,-3309){\makebox(0,0)[lb]{\smash{{\SetFigFont{5}{6.0}{\rmdefault}{\mddefault}{\updefault}{\color[rgb]{0,0,0}$-\eps$}%
}}}}
\put(4166,-2349){\makebox(0,0)[lb]{\smash{{\SetFigFont{5}{6.0}{\rmdefault}{\mddefault}{\updefault}{\color[rgb]{0,0,0}$x$}%
}}}}
\put(9416,-2229){\makebox(0,0)[lb]{\smash{{\SetFigFont{5}{6.0}{\rmdefault}{\mddefault}{\updefault}{\color[rgb]{0,0,0}$y$}%
}}}}
\put(3146,-7284){\makebox(0,0)[lb]{\smash{{\SetFigFont{5}{6.0}{\rmdefault}{\mddefault}{\updefault}{\color[rgb]{0,0,0}$-\eps$}%
}}}}
\put(5036,-7314){\makebox(0,0)[lb]{\smash{{\SetFigFont{5}{6.0}{\rmdefault}{\mddefault}{\updefault}{\color[rgb]{0,0,0}$\eps$}%
}}}}
\put(8396,-7299){\makebox(0,0)[lb]{\smash{{\SetFigFont{5}{6.0}{\rmdefault}{\mddefault}{\updefault}{\color[rgb]{0,0,0}$\frac{1}{2}-\eps$}%
}}}}
\end{picture}%